\documentclass[12pt]{amsart}
\usepackage[colorlinks=true,citecolor=magenta,linkcolor=magenta]{hyperref}
\usepackage{amssymb,verbatim,amscd,amsmath,graphicx}
\usepackage{enumerate}
\usepackage{graphicx}
\usepackage{caption}
\usepackage{subcaption}
\usepackage{pdfsync}
\usepackage{color,colortbl}
\usepackage{fullpage}
\usepackage{bbm}
\usepackage{comment}
\usepackage{tikz}
\numberwithin{equation}{section}
\newtheorem{theorem}{Theorem}[section]

\newtheorem{proposition}[theorem]{Proposition}
\newtheorem{corollary}[theorem]{Corollary}

\theoremstyle{definition}

\newtheorem{conjecture}[theorem]{Conjecture}
\newtheorem{def-prop}[theorem]{Definition-Proposition}

\newtheorem{remark}[theorem]{Remark}
\newtheorem{example}[theorem]{Example}

\newtheorem*{acknowledgement}{Acknowledgements}

\DeclareMathOperator{\reg}{reg}
\DeclareMathOperator{\depth}{depth}

\DeclareMathOperator{\Ass}{Ass}

\newcommand{\CC}{{\mathbb C}}

\newcommand{\PP}{{\mathbb P}}

\newcommand{\NN}{{\mathbb N}}

\def\mm{{\mathfrak m}}

\def\pp{{\frak p}}

\def\ahat{\widehat{\alpha}}

\def\1{{\bf 1}}
\def\0{{\bf 0}}

\begin{document}
	
\title{The Initial Degree of Symbolic Powers of Fermat-like Ideals of Planes and Lines Arrangements}

\author{Th\'ai Th\`anh Nguy$\tilde{\text{\^e}}$n}
\address{Tulane University \\ Department of Mathematics \\
	6823 St. Charles Ave. \\ New Orleans, LA 70118, USA and
	Hue University, College of Education, Vietnam.}
\email{tnguyen11@tulane.edu}
\urladdr{https://sites.google.com/view/thainguyenmath}

\keywords{Fermat Ideals, Fermat Lines Configuration, Resurgence Number, Waldschmidt Constant, Ideals of Points, Ideals of Lines, Symbolic Powers, Containment problem, Stable Harbourne--Huneke Conjecture, Interpolation Problem}
\subjclass[2010]{14N20, 13F20, 14C20}

\begin{abstract}
We explicitly compute the least degree of generators of all symbolic powers of the defining ideal of Fermat-like configuration of lines in $\mathbb{P}^3_\mathbb{C}$, except for the second symbolic powers, where we provide bounds for them. We will also explicitly compute those numbers for ideal determining the singular locus of the arrangement of lines given by the pseudoreflection group $A_3$. As direct applications, we verify Chudnovsky's(-like) Conjecture, Demailly's(-like) Conjecture and Harbourne-Huneke Containment problem as well as calculate the Waldschmidt constant and (asymptotic) resurgence number.
\end{abstract}

\maketitle


\section{Introduction} \label{sec.intro}

Let $n\ge 3$ and $I_n$ be the ideal in $\CC[x,y,z,w]$ generated by: 
$$g_1=(x^n-y^n)(z^n-w^n)xy, \quad g_2=(x^n-y^n)(z^n-w^n)zw $$
$$g_3=(x^n-z^n)(y^n-w^n)xz, \quad g_4=(x^n-z^n)(y^n-w^n)yw $$
$$g_5=(x^n-w^n)(y^n-z^n)xw, \quad g_6=(x^n-w^n)(y^n-z^n)yz $$
which we will refer as \textit{Fermat-like ideal}. This ideal corresponds to the restricted Fermat arrangement of planes in $\PP^3$, where the correspondent variety is the union of all lines with multiplicity at least $3$, i.e, there are at least $3$ planes passing through each line. \textit{Fermat-like ideals} were first introduced by Malara and Szpond in their work \cite{Malara2017FermattypeCO} in an effort to provide counterexamples in higher dimension to the famous containment $I^{(3)} \subseteq I^2$. It can be seen as an analog of \textit{Fermat ideals} in higher dimension. For more information about \textit{Fermat ideals}, we refer interested readers to \cite{NagelSeceleanu}, \cite{counterexamples}, \cite{HaSeFermat}, \cite{JustynaFermat},\cite{ThaiFermat}.\par
\vspace{0.5em}
In this manuscript, we will discover many similarities between these ideals and \textit{Fermat ideals} in terms of the least degree of generators (or initial degree) of their symbolic powers, denoted by $\alpha(I_n^{(k)})$, \textit{Waldschmidt constants} $\ahat(I_n)$, and expectedly, \textit{(asymptotic) resurgence numbers}. Moreover, the verification to certain \textit{Harbourne-Huneke containment} of them can also be decided purely by the knowledge of the aforementioned numerical invariants. This work is a continuation of our paper \cite{ThaiFermat} where we answer those questions for Fermat ideals. Our results on Fermat-like ideals are the following:

\begin{theorem}[Propositions \ref{prop.WaldLineFermat}, \ref{prop.n=2}, Theorems \ref{thm.Fermatlinege5}, \ref{thm.Fermatlinen=4}, \ref{thm.Fermatlinen=3}]
Let $n \ge 3$ be an integer and Fermat-like ideal $I_n$ in $\CC[x,y,z,w]$ as described above. Then:
\begin{enumerate}
    \item For $n \ge 5$, $\alpha(I_n^{(m)})=2nm$ for all $m\ge 3$.
    \item $\alpha(I_4^{(m)})=8m$ for all $m\ge 3$ and $m\not= 5$.
    \item $\alpha(I_3^{(3k)})=18k$ for $k\ge 1$, $\alpha(I_3^{(3k+1)})=18k+8$ and $\alpha(I_3^{(3k+2)})=18k+16$ for $k\ge 0$.
    \item $4n+2\le \alpha(I_n^{(2)})\le 4n+4=\alpha(I_n^{2})$.
    \item $\widehat{\alpha}(I_n) = 2n$.
\end{enumerate}
\end{theorem}
This gives the almost complete answer to the question of the least degree of generators of almost all symbolic powers of the defining ideal of Fermat-like configuration of lines in $\PP^3_\CC$ and bounds only for the case when $m=2$. Note that for each $m,n$ the least generating degree of symbolic power of the Fermat-like ideal is twice as large as that of Fermat ideal in almost all cases (and are expected to be in all cases as suggested by Macaulay2 computation \cite{M2}).\par
\vspace{0.5em}
One important problem that has attracted a lot of attention recently is the \textit{containment problem}, namely, to determine the set of pairs $(m,r)$ for which $I^{(m)}\subseteq I^r$ for a given ideal $I$. Following the celebrated results in \cite{ELS, comparison, MaSchwede} that $I^{(m)}\subseteq I^r$ whenever $m\ge hr$, where $h$ is the \emph{big height} of $I$, the \textit{resurgence} $\rho(I)$ is introduced in \cite{BoH}, and the \textit{asymptotic resurgence} $\widehat{\rho}(I)$ is introduced in \cite{AsymptoticResurgence} in order to study the pairs $(r, m)$ numerically and has turned out to be very useful invariants.
$$\rho(I):= \sup \{ \dfrac{m}{r} : I^{(m)} \not\subseteq I^r \}$$
$$\widehat{\rho}(I):= \sup \{ \dfrac{m}{r} : I^{(mt)} \not\subseteq I^{rt} \text{  for  } t\gg 0 \}.$$

It is well-known that over a regular ring, $1 \leqslant \rho(I) \leqslant h$, where $h$ denotes the big height of $I$, that is the maximum height of its associated primes. As noted in \cite[Remark 2.7]{GrifoStable}, the Stable Harbourne Conjecture, namely, $I^{(hr-h+1)} \subseteq I^r$ for $r\gg 0$, follows immediately whenever $\rho(I) < h$. In that case, we say that $I$ has \emph{expected resurgence} \cite{GrifoHunekeMukundan}. We show that the defining ideals of the singular locus of the line arrangements corresponds to the group $A_3$ and $B_3$ have expected resurgence and strongly predict that the same applies to $I_n$. \par
\vspace{0.5em}
Another motivation for computing the least degree of generators of their symbolic powers is to provide more evidence for ideals that satisfy \textit{Chudnovsky-like inequality} and \textit{Demailly-like inequality} as well as \textit{Harbourne-Huneke Containment} and stable containment. This motivation stems from the previous work by Bisui, Grifo, Hà and the author, see \cite[Section 3]{BGHN2-2020} and containment in \cite{BGHN2020}. Similar to the Fermat ideals, the verification and failure of the containment for $I_n$ can be checked purely by their numerical invariants including the least degree of generators of symbolic powers, the regularity, or the resurgence number and the maximal degree of generators. As direct applications, we will provide some evidence for classes of ideals beyond points that satisfy certain Harbourne-Huneke Containment (a similar question appears as $3.1$ in \cite[Section 3]{BGHN2-2020}, that all Fermat-like ideals satisfy the \textit{Harbourne-Huneke Containment}. Note that the big height of the ideal $I_n$ is 2 for each $n$, see the description of $I_n$ in section 3.

\begin{corollary}[Corollary \ref{cor.HaHu}]
The Fermat-like ideals verify the following Harbourne-Huneke containment (see \cite[Conjecture 2.1]{HaHu} for ideal of points) $$I_n^{(2r)} \subseteq \mm^rI_n^r, \quad \forall r.$$
\end{corollary}

As another direct application of these computations, we will give affirmed examples for question $3.2$ in \cite[Section 3]{BGHN2-2020}, that all Fermat-like ideals satisfy the \textit{Chudnovsky-like inequality} and \textit{Demailly-like inequality}.

\begin{corollary}
Fermat-like ideals satisfy Demailly-like inequality (and hence, Chudnovsky-like inequality) $$\widehat{\alpha}(I_n)\ge \frac{\alpha(I_n^{(m)}) + h - 1}{m + h -1}$$ for all $m$, where $h=2$.
\end{corollary}

Along the way, we give the free resolution, Castelnuovo-Mumford regularity formula for $I_n$ as a step towards investigating Harbourne-Huneke containment and computing (asymptotic) resurgence number for $I_n$.\par
\vspace{0.5em}
In previous work \cite{ThaiFermat}, we investigate the above questions for Fermat ideals. In particular, the Fermat ideal for $7$ points can be seen as ideal of the singular locus of the arrangement of lines given by the pseudoreflection group $G(2,2,3)=D_3$. We continue to study the defining ideal of the singular locus of the arrangement of lines given by the group $G(1,1,4)=A_3$ in this paper and give some discussion for group $G(2,1,3)=B_3$ (see \cite{DrabkinSeceleanu} for description of these groups and related line arrangements). These are $3$ groups that have correspondent ideals with small degree of generators, that had to be considered separately in the result involving Harbourne containment in \cite[Proposition 6.3]{DrabkinSeceleanu}. It turns out interestingly that the ideal corresponding to the group $A_3$ is the same as that of group $D_3$ in term of least degree of generators of symbolic powers, Waldschmidt constant, (asymptotic) resurgence numbers; and hence, satisfies Chudnovsky's Conjecture and Demailly's Conjecture as well as Harbourne-Huneke Containment, stable Harbourne Containment and some stronger containment.

\begin{theorem}\label{thm.idealJ}
Let $J$ be the ideal of the singular locus of the arrangement of lines given by the group $G(1,1,4)=A_3$. Then 
\begin{enumerate}
    \item $\alpha(J^{(2k)})=5k$ for all $k\ge 2$ and $\alpha(J^{(2)})=6$.
    \item $\alpha(J^{(2k+1)})=5k+3$ for all $k\ge 0$.
    \item $\widehat{\alpha}(J) = \dfrac{5}{2}$.
    \item $\rho(J)=\widehat{\rho}(J)=\dfrac{6}{5}$.
\end{enumerate}
\end{theorem}

Similar to the Fermat ideals and Fermat-like ideals, the verification and failure of the containment $J$ can be checked purely by their numerical invariants including the least degree of generators of symbolic powers, the regularity, or the resurgence number and the maximal degree of generators.

\begin{corollary}
The ideal $J$ in theorem \ref{thm.idealJ} verifies the following containment and stable containment
\begin{enumerate}
    \item Harbourne-Huneke containment (see \cite[Conjecture 2.1]{HaHu}) $$J^{(2r)} \subseteq \mm^rJ^r, \quad \forall r$$
    \item Harbourne-Huneke containment (see \cite[Conjecture 4.1.5]{HaHu}) $$J^{(2r-1)} \subseteq \mm^{r-1}J^r, \quad \forall r \ge 1$$ 
    \item A stronger containment $J^{(2r-2)} \subseteq \mm^{r}J^r, \quad \forall r\ge 5$
\end{enumerate}
\end{corollary}
 
We work over the field of complex numbers but our results hold over any algebraically closed field of characteristic $0$.

\begin{acknowledgement}
	The author would like to thank his advisor, Tài Huy Hà, for introducing him this subject and giving many helpful suggestions and comments. He also thanks Alexandra Seceleanu, Abu C. Thomas, Benjamin Drabkin, and Paolo Mantero for comments on an early draft of this manuscript. Finally, he thanks the referee for a careful reading of the paper and for many valuable suggestions.
\end{acknowledgement}


\section{Preliminaries}

Let $R = \CC[x_0,\ldots ,x_N]$ be the homogeneous coordinate ring of $\PP^N$, and let $\mm$ be its maximal homogeneous ideal. For a homogeneous ideal $I \subseteq R$, let $\alpha(I)$ denote the least degree of a non-zero homogeneous polynomial in $I$, and let 
$$I^{(m)} := \bigcap_{\pp \in \Ass(R/I)} I^m R_\pp \cap R$$
denote its $m$-th \emph{symbolic power}.\par
\vspace{0.5em}
The \textit{Waldschmidt constant} of $I$ is defined to be the limit and turned out to be the infimum 
$$\widehat{\alpha}(I):= \lim_{m\rightarrow \infty} \dfrac{\alpha(I^{(m)})}{m} = \inf _{m\rightarrow \infty} \dfrac{\alpha(I^{(m)})}{m}$$ 

In studying ideals defining sets of points in $\PP^N$, \textit{Chudnovsky's Conjecture} gives a lower bound for $\alpha(I^{(m)})$ in term of $\alpha(I),N$ and $m$ as follows.

\begin{conjecture}[Chudnovsky]
	\label{conj.Chud}
	Let $I$ be the defining ideal of a set of points in $\PP^N_\CC$. Then, for all $n \ge 1$,
$$\dfrac{\alpha(I^{(n)})}{n} \ge \dfrac{\alpha(I)+N-1}{N}$$
\end{conjecture}

It is also natural to ask if the \textit{Chudnovsky-like inequality} is still true for a homogeneous radical ideal $I$, for example, if we replace $N$ by big height $h$ of $I$. Many partial results are known for Chudnovsky's Conjecture, for example, in \cite{EsnaultViehweg, BoH, HaHu, GHM2013, Dumnicki2015, DTG2017, FMX2018, BGHN2020}. Recently, the conjecture was proved for a \emph{very general} set of points in \cite{DTG2017, FMX2018}, for a \emph{general} set of sufficiently many points in \cite{BGHN2020}, and for small numbers of general points in \cite{SankhoThaiChudnovsky}. On the other hand, Chudnovsky-like inequality was verified for a set of \textit{very general} lines in $\PP^3$ in \cite{ChudnovskyGeneralLines}. There is also result that for any homogeneous ideal $I$, ideal $I^{(t)}$ satisfies Chudnovsky-like inequality for all $t\gg 0$, see \cite{FMX2018}. A defining ideal of fat points with equi-multiplicity at least $2$ is also shown to satisfy Chudnovsky-like inequality, see \cite{BGHN2020}. The following generalization is due to Demailly \cite{Demailly1982}.

\begin{conjecture}[Demailly] \label{conj.Demailly}
	Let $I$ be the defining ideal of a set of points in $\PP^N_\CC$. Let $m \in \NN$ be any fixed integer. Then, for all $n \ge 1$,
$$\dfrac{\alpha(I^{(n)})}{n} \ge \dfrac{\alpha(I^{(m)}) + N-1}{m+N-1}$$ 
\end{conjecture}

\textit{Demailly's Conjecture} for $N = 2$ was proved by Esnault and Viehweg \cite{EsnaultViehweg}. In higher dimension, Demailly's Conjecture holds for a \emph{very general} set of sufficiently many points (the number of points depends on $m$) in \cite{MSS2018} and for a \emph{general} set of $k^N$ points in $\PP^N$ \cite{CJ2020}. In \cite{BGHN2-2020}, the conjecture was proved for a \emph{general} set of sufficiently many points (also depends on $m$). We also raise a question \cite[Question 3.2]{BGHN2-2020} for ideal that verified Demailly-like inequality and show some examples of such ideals including: defining ideal of a codimension $h$ star configuration in $\PP^N$, generic determinantal ideals, determinantal ideals of symmetric matrices and pfaffian ideals of skew symmetric matrices.\par
\vspace{0.5em}
The \textit{resurgence number} $\rho(I)$ is introduced in \cite{BoH} as
$$\rho(I):= \sup \{ \dfrac{m}{r} : I^{(m)} \not\subseteq I^r \}$$
and the \textit{asymptotic resurgence} $\widehat{\rho}(I)$ is introduced in \cite{AsymptoticResurgence} as
$$\widehat{\rho}(I):= \sup \{ \dfrac{m}{r} : I^{(mt)} \not\subseteq I^{rt} \text{  for  } t\gg 0 \}$$
in effort to study the \textit{containment problem} numerically, namely, the pairs $(m,r)$ such that $I^{(m)}\subseteq I^r$. It is a celebrated result in \cite{ELS, comparison, MaSchwede} that $I^{(hr)}\subseteq I^r$ for codimension $h$ ideal $I$. In another effort to improve this containment for ideal of points and study Chudnovsky's Conjecture and Demailly's Conjecture, Harbourne and Huneke in \cite{HaHu} conjectured that the defining ideal $I$ for any set of points in $\PP^N$ satisfies some stronger containment $$I^{(Nm)} \subseteq \mm^{m(N-1)}I^m,  \text{  and   } I^{(Nm-N+1)} \subseteq \mm^{(m-1)(N-1)}I^m,$$
for all $m \geqslant 1$. We also raise a question \cite[Question 3.1]{BGHN2-2020} for other classes of ideals that satisfies some version of \textit{Harbourne-Huneke Containment} (replacing $N$ by $h$) and showed some examples including the ideals mentioned in previous paragraph. Note that suitable version of \textit{stable Harbourne-Huneke Containment} (or even infinitely many such containment) would imply Chudnovsky's and Demailly' Conjecture (or Demailly-like inequality), for example, see \cite{BGHN2-2020}. On the other hand, the resurgence always satisfies $\rho(I) \geqslant 1$, and over a regular ring, the resurgence of a radical ideal is always at most the big height $h$. As noted in \cite[Remark 2.7]{GrifoStable}, the Stable Harbourne Conjecture $I^{(hr-h+1)} \subseteq I^r$ for $r\gg 0$ follows immediately whenever $\rho(I) < h$. In that case, we say that $I$ has \emph{expected resurgence} \cite{GrifoHunekeMukundan}.We refer interested readers to \cite{CHHVT2020} for more information about the Waldschmidt constant, resurgence number, containment between symbolic and ordinary powers of ideals.

\section{Ideal of restricted Fermat configuration of lines in $\PP^3$}
Recall that Fermat-like ideal $I_n$ is the ideal generated by 
$$g_1=(x^n-y^n)(z^n-w^n)xy, \quad g_2=(x^n-y^n)(z^n-w^n)zw $$
$$g_3=(x^n-z^n)(y^n-w^n)xz, \quad g_4=(x^n-z^n)(y^n-w^n)yw $$
$$g_5=(x^n-w^n)(y^n-z^n)xw, \quad g_6=(x^n-w^n)(y^n-z^n)yz $$

From \cite{Malara2017FermattypeCO}, geometrically, $I_n$ is the defining ideal of the union of lines with multiplicity at least $3$ of the Fermat arrangement of flats (planes) in $\PP^3$, denoted by $H_j$, $j=1,\ldots,6n$, that is defined by the vanishing of the polynomial: $$F_n=(x^n-y^n)(z^n-w^n)(x^n-z^n)(y^n-w^n)(x^n-w^n)(y^n-z^n)$$ 
There are $4n^2+6$ lines, denoted by $L_j$, $j=1,\ldots,4n^2+6$, in the above restricted Fermat configuration of lines. Notice also that from this description, each of $6n$ planes $H_j$ passes through exactly $2n+1$ lines, for example, the plane $x=\epsilon y$ passes through the line defined by $(x,y)$, $n$ lines defined by $(x-\epsilon y,x- \epsilon ^k z)$ and $n$ lines defined by $(x-\epsilon y,y- \epsilon ^k w)$ for $k=0,\ldots n-1$, where $\epsilon$ is an $n-$root of $1$.\par
\vspace{0.5em}
Algebraically, let $f_n=(x^n-y^n)(z^n-w^n), g_n=(x^n-z^n)(y^n-w^n), h_n=(x^n-w^n)(y^n-z^n)$ and $K_n=(f_n,g_n)$, we can write 
$$I_n=K_n \cap (x,y) \cap (x,z) \cap (x,w) \cap (y,z) \cap (y,w) \cap (z,w)$$ 
and since $f_n,g_n$ form a regular sequence, for any $m\ge 1$ we have
$$I_n^{(m)}=K_n^m \cap (x,y)^m \cap (x,z)^m \cap (x,w)^m \cap (y,z)^m \cap (y,w)^m \cap (z,w)^m$$

In particular, for each $n$, the ideal $I_n$ is a radical ideal with big height $2$. In the following, we will compute the initial degree of symbolic powers of $I_n$. Recall our main strategy in previous work \cite{ThaiFermat}, we will study a subsequence of $\alpha(I_n^{(m)})$, which gives us information about $\widehat{\alpha}(I_n)$, then use this to calculate other $\alpha(I_n^{(m)})$.

\begin{proposition}\label{prop.WaldLineFermat}
For $n\ge 3$, we have $\widehat{\alpha}(I_n) = 2n$.
\end{proposition}
\begin{proof}
For every $m\ge 1$, since $I_n^{(m)} \subset K_n^m$, we have $\alpha(I_n^{(m)})\ge \alpha(K_n^m)=2nm$. Thus, $\widehat{\alpha}(I_n) \ge 2n$. On the other hand, since $h_n=g_n-f_n \in K_n$, we have $F_n^m =f_n^m g_n^m h_n^m \in K_n^{3m}$. Furthermore, since $nm\ge 3m$, $F_n^m \in (x,y)^{3m} \cap (x,z)^{3m} \cap (x,w)^{3m} \cap (y,z)^{3m} \cap (y,w)^{3m} \cap (z,w)^{3m}$ so $F_n^m \in I_n^{(3m)}$,  hence $\alpha(I_n^{(3m)}) \le 6nm $ for every $m\ge 1$. Therefore $\widehat{\alpha}(I_n)\le 2n$ and it follows that $\widehat{\alpha}(I_n) = 2n$.
\end{proof}

Now following exactly the same argument as those of \cite[Section 3]{ThaiFermat} we have the following two results:

\begin{theorem}\label{thm.Fermatlinege5}
For $n\ge 5$, we have $$\alpha(I_n^{(m)})=2nm$$ for all $m\ge 3$.
\end{theorem}

\begin{proof}
The proof is identical to that of \cite[Theorem 3.1]{ThaiFermat}. The idea is that we have $f_ng_nh_n(f_n,g_n)^{m-3} \subseteq I_n^{(m)}$ for $3\le m \le n$ and that $(f_ng_nh_n)^k(f_n,g_n)^{k(n-3)-a} \subseteq I^{(kn-a)}$ for $k\ge 2$, $0\leq a \leq n-1$.
\end{proof}

\begin{theorem}\label{thm.Fermatlinen=4}
For $n=4$, for all $m\ge 3$ but $m \not = 5$, $$\alpha(I_4^{(m)})=8m$$
\end{theorem}

\begin{proof}
The proof is identical to that of \cite[Theorem 3.2]{ThaiFermat}. The idea is the same as that of above theorem. For $m=5$, we can check by the same argument that the element $yzf_4g_4^2h_4^2 \in I_4^{(5)}$. Macaulay2 computation \cite{M2} suggests that, in fact, $\alpha(I_4^{(5)})=42$.
\end{proof}

Now we calculate the least degree of generators of symbolic powers for $I_3$.

\begin{theorem}\label{thm.Fermatlinen=3}
For $n=3$ and for $m\ge 1$ we have the following
\begin{enumerate}
    \item $\alpha(I_3^{(3m)})=18m$ for $m\ge 1$,
    \item $\alpha(I_3^{(3m+1)})=18m+8$ for $m\ge 0$,
    \item $\alpha(I_3^{(3m+2)})=18m+16$ for $m\ge 0$.
\end{enumerate}
\end{theorem}

\begin{proof}
\begin{enumerate}
    \item Since $\widehat{\alpha}(I_3)=6$, $\alpha(I_3^{(3m)})\ge 18m$ for $m\ge 1$. As we saw earlier, $F_3^m \in I_3^{(3m)}$ for $m\ge 1$ so $\alpha(I_3^{(3m)})\le 18m$ and thus $\alpha(I_3^{(3m)}) = 18m$ for $m\ge 1$.
    \item Suppose that there is $m\ge 1$ such that $\alpha(I_3^{(3m+1)})\le 18m+7$. Then there is a divisor $D$ of degree $18m+7$ vanishing to order at least $3m+1$ along every line of $42$ lines $L_i$ in the restricted Fermat configuration. Intersecting $D$ with any of the $18$ planes $H_j$, $j=1,\ldots,18$, since each planes $H_j$ contains exactly $7$ lines, suppose that $D$ doesn't contain $H_j$ then by the generalized Bezout Theorem, the intersection of $D$ and $H_j$ is of dimension $1$ and degree $18m+7$, which is a contradiction since $18m+7<7(3m+1)$. \par
    \vspace{0.5em}
    Thus we conclude that each $H_j$ is a component of $D$ for all $j=1,\ldots,18$. Hence, there exists a divisor $D'=D- \sum_{j=1}^{18} H_j$ of degree $18(m-1)+7$ vanishing to order at least $3(m-1)+1$ along every line $L_i$. Repeating this argument $m$ times we get a contradiction with $\alpha(I_3)=8$. Thus $\alpha(I_3^{(3m+1)})\ge 18m+8$ for all $m \ge 1$. On the other hand, it easy to see that, $f_3^mg_3^mh_3^{m+1}yz\in I_3^{(3m+1)}$ so $\alpha(I_3^{(3m+1)})\le 18m+8$ for all $m \ge 1$. Therefore, $\alpha(I_3^{(3m+1)})= 18m+8$ for all $m \ge 1$.
    \item By identical argument to previous case and the fact that $\alpha(I_3^{(2)})=16$, we have that $\alpha(I_3^{(3m+2)})\ge 18m+16$ for all $m \ge 1$. Since the element $f_3^{m+1}g_3^{m+1}h_3^mxyzw \in I_3^{(3m+2)}$ we have that $\alpha(I^{(3m+2)}) \le 18m+16$ for all $m\ge 1$.
\end{enumerate}

\end{proof}

From above computation, we showed that for each $m,n$, the least generating degree of symbolic power of the Fermat-like ideal is twice as large as that of Fermat ideal \cite{ThaiFermat} in almost all cases (we have not shown for $m=2$, but it is expected to be the case as well).

\begin{proposition}\label{prop.n=2}
For all $n\ge 3$, $4n+2 \le \alpha(I_n^{(2)})\le 4n+4$.
\end{proposition}

\begin{proof}
We know that $\alpha(I_n^{(2)}) \le \alpha(I_n^2) =2(2n+2)$. Now suppose that $\alpha(I_n^{(2)}) \le 4n+1$. Then there is a divisor $D$ of degree $4n+1$ vanishing to order at least $2$ along every line $L_i$ in the restricted Fermat configuration. Since the intersection of $D$ and any plane $H_j$ consists of $2n+1$ lines to order at least $2$, by generalized Bezout theorem, each $H_j$ is a component of $D$ because $\deg(D).\deg(H_j)=4n+1<2(2n+1)$. This is a contradiction since there are $6n$ planes and $6n>4n+1=\deg(D)$ when $n\ge 3$. Therefore, $4n+2 \le \alpha(I_n^{(2)})\le 4n+4$ for all $n\geq 3$.
\end{proof}

\begin{remark}
Macaulay2 computations \cite{M2} for $n$ small suggest that $\alpha(I^{(2)}) =2(2n+2)$ for $n\ge 3$.
\end{remark}

As mentioned earlier, in \cite[Question 3.2]{BGHN2-2020}, we raise a question to study the Demailly-like bound for homogeneous ideals. We also give example of classes of ideals that satisfy Demailly-like inequality. From the above calculations, as a direct consequence, we see that the ideal $I_n$ also gives an example for Demailly-like inequality. 

\begin{corollary}
For $n\ge 3$, the ideal $I_n$ satisfies Demailly-like inequality $$\widehat{\alpha}(I_n) \ge \dfrac{\alpha(I_n^{(m)})+h-1}{m+h-1}$$ for all $m\ge 1$ where $h$ is the big height of $I_n$.
\end{corollary}
\begin{proof}
Direct from the above calculation with notice that $h=2$.
\end{proof}

The following containment are also direct consequences of the above calculations.

\begin{corollary}\label{cor.HaHu}
For every $n\geqslant 3$, restricted Fermat configuration ideal verifies Harbourne-Huneke containment $$I_n^{(2r)} \subseteq \mm^rI_n^r, \quad \forall r$$
\end{corollary}

\begin{proof}
Since $I_n$ is a radical ideal with big height $2$, by \cite{ELS}, we have $I_n^{(2r)} \subseteq I_n^r, \forall r$, so the above containment come from the fact that $$\alpha(I_n^{(2r)})\geqslant r+\omega(I_n^r)$$ for all $n\geqslant 3$ and $r \geqslant 1$ (the case $r=0$ is trivial). Indeed, $\alpha(I_n^{(2r)})\geqslant 4rn >r+r(2n+2) = r+\omega(I_n^r)$ for $n\geqslant 3$.
\end{proof}

In the works of studying containment for ideals of points, or more general, for ideals of smooth schemes, the result by Bocci and Harbourne (\cite[Lemma 2.3.4]{BoH}), which says that if $\alpha(I^{(m)}) \geqslant \reg(I^r)$ then $I^{(m)} \subseteq I^r$, turns out to be extremely useful (see for example, \cite{BoH,BGHN2020,BGHN2-2020}. Now we calculate the regularity of the ideal $I_n$ for all $n\geqslant 3$. 

\begin{proposition}
For all $n\geqslant 3$, $R/I_n$ is Cohen-Macaulay and its minimal graded free resolution is
$$0 \to \begin{matrix} R(-(2n+3))^{4} \\
\oplus \\ R(-4n)\end{matrix} \stackrel{\varphi_2}\to R(-(2n+2))^6 \stackrel{\varphi_1} \to R \to R/I_n\to 0 $$
where $\varphi_1 = [g_1,g_2,g_3,g_4,g_5,g_6]$ and 
$$\varphi_2 = \begin{bmatrix}
-y&0&-w&0&x^{n-1}z^{n-1} \\
x &0 &0&-w&-y^{n-1}z^{n-1} \\
z&-w&0&0&-x^{n-1}y^{n-1} \\ 
0&-y&z&0&-x^{n-1}w^{n-1} \\ 
0&x&0&z&y^{n-1}w^{n-1} \\ 
0&0&x&-y&-z^{n-1}w^{n-1}
\end{bmatrix}$$
The maximal minors of $\varphi_2$ are multiple of $g_1,\ldots , g_6$ as in Hilbert-Burch theorem.
\end{proposition}

\begin{proof}
Direct calculation give the minors of $\varphi_2$ deleting row $1,2,3,4,5,6$ are $2g_3$, $-2g_6$, $2g_1$, $-2g_5$, $2g_4$, $-2g_2$ respectively. Applying \cite[Theorem 18.18]{Eisenbud}, we see that $I_n$ is the ideal generated by $5\times 5$ minors of the $6\times 5$ matrix $\varphi_2$ and since $I_n$ has codimension $2$ which agrees with $(6-5+1)(5-5+1)$, $R/I_n$ is Cohen-Macaulay. In particular, $\depth(R/I_n)=2$ and by  Auslander–Buchsbaum formula $R/I_n$ has projective dimension $2$. Thus by Hilbert-Burch theorem, $R/I_n$ has such minimal graded free resolution.
\end{proof}

The immediate consequence of the above free resolution is the regularity of the ideal $I_n$.

\begin{corollary}
For all $n\geqslant 3$, $\reg(I_n)=4n-1$. 
\end{corollary}

\begin{proposition}\label{resurgenceI3}
For all $m,r$ such that $\dfrac{m}{r}> \dfrac{3}{2}$, we have the inequality $\alpha(I_3^{(m)}) \geqslant \reg(I_3^r)$.
\end{proposition}

\begin{proof}
We will apply \cite[Theorem 0.6]{chardin2007behavior}, which states the following:\par
\vspace{0.5em}
Let $I$ be an homogeneous ideal of a polynomial ring $R$ over a ﬁeld, generated in degrees at most $d$, such that $\dim(R/I) = 2$. Assume that $I_p \subseteq R_p$ is a complete intersection for every prime $I \subseteq p$ such that $\dim R/p = 2$. Then
$$\reg(I^2) \leqslant \max \{ 2\reg(I), \reg(I^{sat})+2d-2 \}$$
and for $r\geqslant 3$
$$\reg(I^r) \leqslant \max \{ 3\reg(I)+(r-3)d, \reg(I^{sat})+rd-2 \}.$$ 
Clearly, $\dim(R/I_3)=2$. Since the set of associated primes of $I_3$ is also the set of its minimal primes which consists of $42$ minimal primes $P_j$ where each minimal prime is the defining ideal of one of a line in the configuration. Since any prime ideal $P$ that contains $I_3$ must contain one $P_j$, such prime ideal $P$ with $\dim(R/P)=2$ has to be $P_j$. Thus, for any prime $I\subset P$ such that $\dim(R/P)=2$, $(I_3)_{P}$ is a complete intersection. Since $\reg(I_3)=11$, $\alpha(I_3)=8$ and $I_3=I_3^{sat}$, by \cite[Theorem 0.6]{chardin2007behavior} we have for $r\geqslant 3$:
$$\reg(I_3^r) \leqslant \max\{ 3\reg(I_3)+8(r-3),\reg(I_3)+8r-2 \}=\max \{ 33+8(r-3),11+8r-2 \}=8r+9$$
and $\reg(I_3^2) \leqslant \max \{22,11+16-2  \}=25$. Thus, for every $r\geqslant 2$: $$\reg(I_3^r) \leqslant 8r+9$$
Now for $\dfrac{m}{r}> \dfrac{3}{2}$, i.e, $2m \geqslant 3r+1$, we have $\alpha(I_3^{(m)}) \geqslant 6m \geqslant 9r+3$. Since $$9r+3 \geqslant 8r +9 \geqslant \reg(I_3^r) $$ for all $r\geqslant 6$, it suffices now to check that for $r\leqslant 5$ and $2m \geqslant 3r+1$ we still have $\alpha(I_3^{(m)}) \geqslant \reg(I_3^r)$. Indeed,
\begin{itemize}
    \item When $r=1$, for $m\geqslant 2$, $\reg(I_3)=11 < 16 = \alpha(I_3^{(2)}) \leqslant \alpha(I_3^{(m)})$.
    \item When $r=2$, for $m\geqslant 4$, $\reg(I_3^2)\leqslant 25 < 26 = \alpha(I_3^{(4)}) \leqslant \alpha(I_3^{(m)})$.
    \item When $r=3$, for $m\geqslant 5$, $\reg(I_3^3)\leqslant 33 < 34 = \alpha(I_3^{(5)}) \leqslant \alpha(I_3^{(m)})$.
    \item When $r=4$, for $m\geqslant 7$, $\reg(I_3^4)\leqslant 41 < 44 = \alpha(I_3^{(7)}) \leqslant \alpha(I_3^{(m)})$.
    \item When $r=5$, for $m\geqslant 8$, $\reg(I_3^5)\leqslant 49 < 52 = \alpha(I_3^{(8)}) \leqslant \alpha(I_3^{(m)})$.
\end{itemize}
\end{proof}

\begin{corollary}
Suppose that the inequality $\alpha(I_n^{(m)})\geqslant \reg(I_n^r)$ implies that $I_n^{(m)} \subseteq I_n^r$. Then $\rho(I_3)=\dfrac{3}{2}$.
\end{corollary}
\begin{proof}
By \cite{Malara2017FermattypeCO}, $I_n^{(3)} \not \subseteq I_n^2$ hence, $\rho(I_3) \geqslant \dfrac{3}{2}$. By Proposition \ref{resurgenceI3}, for all $m,r$ such that $\dfrac{m}{r}> \dfrac{3}{2}$, we have $I^{(m)} \subseteq I^r$. Hence, $\rho(I_3) \leqslant \dfrac{3}{2}$.
\end{proof}

We also deduce the following results as a step towards calculating the resurgence number of $I_n$ for $n\geqslant 4$.

\begin{corollary}
If $n\geqslant 4$ then for all $m,r$ such that $\dfrac{m}{r}>\dfrac{3}{2}$ and $r\geqslant 6$, we have the inequality $\alpha(I_n^{(m)})\geqslant \reg(I_n^r)$. 
\end{corollary}

\begin{proof}
By the same argument as that of in the proof of proposition \ref{resurgenceI3}, $I_n$ satisfies all conditions of \cite[Theorem 0.6]{chardin2007behavior}, so we have
$$\reg(I_n^r) \leqslant \max \{3(4n-1)+(r-3)(2n+2),4n-1+r(2n+2)-2 \} = r(2n+2)+6n-9 $$
for $r\geqslant 3$. Since $\alpha(I_n^{(m)})\geqslant 2nm$, for $m,r$ such that $\dfrac{m}{r}>\dfrac{3}{2}$, i.e $2m \geqslant 3r+1$ we have $\alpha(I_n^{(m)})\geqslant n(3r+1)$. Now the inequality $$n(3r+1) \geqslant r(2n+2) +6n-9 \Leftrightarrow r\geqslant \dfrac{5n-9}{n-2}=5+\dfrac{1}{n-2},$$
is true for $r\geqslant 6$ and $n\geqslant 4$.
\end{proof}

\begin{remark}
For $n\geqslant 4$, suppose that the inequality $\alpha(I_n^{(m)})\geqslant \reg(I_n^r)$ implies that $I_n^{(m)} \subseteq I_n^r$. Then the resurgence number $\rho(I_n)$ can only be one of the following numbers: $\frac{3}{2},\frac{5}{3},\frac{7}{4}, \frac{8}{5}$ or $\dfrac{9}{5}$. More precisely, among the four containment $I_n^{(9)} \subseteq I_n^5$, $I_n^{(7)} \subseteq I_n^4$, $I_n^{(5)} \subseteq I_n^3$ and $I_n^{(8)} \subseteq I_n^5$, the ratio $\dfrac{m}{r}$ of the first containment (in that order) that fails is the exact value for $\rho(I_n)$, otherwise, $\rho(I_n)=\frac{3}{2}$. \par
\vspace{0.5em}
In fact, we know that $I_n^{(3)} \not \subseteq I_n^2$ by \cite{Malara2017FermattypeCO}, hence, $\rho(I_n) \leqslant \frac{3}{2}$. On the other hand, $I_n^{(m)} \subseteq I_n^r$ for all $m,r$ such that $\frac{m}{r}>\frac{3}{2}$ and $r\geqslant 6$.
The containment are still true if $\frac{m}{r}>\frac{3}{2}$ for $r=1,2$ since $\frac{m}{r}>\frac{3}{2}$ implies that $m\geqslant 2r$ in these cases and $I_n$ has big height $2$. Moreover, if $r=3$ then $\frac{m}{r}>\frac{3}{2}$ implies $m\geqslant 5$ and we know that $I_n^{(m)} \subseteq I_n^3$ for $m\geqslant 6$. Similarly, if $r=4$ then $\frac{m}{r}>\frac{3}{2}$ implies $m\geqslant 7$ and $I_n^{(m)} \subseteq I_n^4$ for $m\geqslant 8$. Lastly, if $r=5$ then $\frac{m}{r}>\frac{3}{2}$ implies $m\geqslant 8$ and $I_n^{(m)} \subseteq I_n^5$ for $m\geqslant 10$.\par
\vspace{0.5em}
If the following four containment $I_n^{(9)} \subseteq I_n^5$, $I_n^{(7)} \subseteq I_n^4$, $I_n^{(5)} \subseteq I_n^3$ and $I_n^{(8)} \subseteq I_n^5$ all hold, then $I_n^{(m)} \subseteq I_n^r$ for all $m,r$ such that $\frac{m}{r}>\frac{3}{2}$ (because $\frac{3}{2}<\frac{8}{5} <\frac{5}{3}<\frac{7}{4}< \frac{9}{5}$), and it follows that $\rho(I_n)=\frac{3}{2}$. Otherwise, suppose $I_n^{(m_0)} \subseteq I_n^{r_0}$ is the first containment that fails among the above four containment (in order), then $\rho(I_n) \geqslant \frac{m_0}{r_0}$. For all $m,r$ such that $\frac{m}{r}>\frac{m_0}{r_0} >\frac{3}{2}$, we know that either $I_n^{(m)} \subseteq I_n^r$ if $\frac{m}{r}$ is not among $\frac{5}{3},\frac{7}{4}, \frac{8}{5}$ or $\frac{9}{5}$ by above result or $\frac{m}{r}$ is among $\frac{5}{3},\frac{7}{4}, \frac{8}{5}$ or $\frac{9}{5}$ that is greater than $\dfrac{m_0}{r_0}$, in which $I_n^{(m)} \subseteq I_n^r$ hold as well. Either way we have $\rho(I_n) \leqslant \frac{m_0}{r_0}$, therefore $\rho(I_n) = \frac{m_0}{r_0}$.\par
\vspace{0.5em}
Macaulay2 calculations \cite{M2} for $n,r$ small suggests that $\reg(I_n^r)=(r+1)\alpha(I_n)-5=(r+1)(2n+2)-5$, and this would imply $\alpha(I_n^{(m)}) \geqslant \reg(I_n^r)$ for all $\frac{m}{r}> \frac{3}{2}$. It is expected that $\rho(I_n)=\dfrac{3}{2}$ for all $n\geqslant 3$.
\end{remark}

It is expected that the least generating degree of symbolic power of the Fermat-like ideal is twice as large as that of the Fermat ideal for each $m,n$ and that they are expected to have the same (asymptotic) resurgence number for each $n$. It would be interesting to see if this is still the case when we generalize the restricted Fermat configuration into higher dimension. It is worth to note that the above calculation suggests that the Fermat-like ideals have expected resurgence. Table \ref{table:Fermat} summarizes all numerical values for the initial degrees of symbolic powers of the ideals $I_n$ that we have calculated in this section.

\begin{table}[h!]
{\small
  \begin{tabular}{|p{1.2cm}||p{0.8cm}|p{1cm}|p{1.2cm}|p{1.4cm}||p{0.8cm}|p{0.8cm}|p{1.2cm}||p{1.8cm}|p{1cm}|}
    \hline
    $n$ &
      \multicolumn{4}{c||}{3}&
      \multicolumn{3}{c||}{4}&
      \multicolumn{2}{c|}{$n \geqslant 5$} \\
      \hline
    $m$ & 2 & $3k$ & $3k+1$ & $3k+2$ & 2 & 5 & $\geqslant 3$, $\not= 5$ & 2 & $\geqslant 3$ \\
    \hline
    $\alpha(I_n^{(m)})$ & 16 & $18k$ & $18k+8$ & $18k+16$ & 20 & 42 & $4m$ & $\geqslant 4n+2$, $\leqslant 4n+4$ & $2nm$ \\
    \hline
    $\widehat{\alpha}(I_n)$ &
      \multicolumn{4}{c||}{3}&
      \multicolumn{3}{c||}{4}&
      \multicolumn{2}{c|}{$n$} \\
    \hline
  \end{tabular}
  \caption{Least degree of generators of symbolic powers of $I_n$}
  \label{table:Fermat}}
\end{table}

\section{Arrangements given by the group $A_3$}
In this section, we will deal with ideal $J=(yz(y-z),zx(z-x),xy(x-y))$, which is the ideal of the singular locus of the arrangement of lines given by the group $G(1,1,4)=A_3$. In general, given a finite group $G \subseteq \text{GL}_{n+1}(\CC)$ generated by pseudoreflections. By a pseudoreflection we mean a linear transformation that is not the identity and fixes a hyperplane pointwise and have finite order. We can regard the generators of $G$ as a hyperplane arrangement where the hyperplanes are pointwise fixed by the elements of $G$ that are pseudoreflections. It is shown in \cite[Proposition 3.8]{DrabkinSeceleanu} that the singular locus (the Jacobian ideal) of the arrangement of lines correspond to $G(1,1,4)=A_3$ is given by $J$. We will see that its least degree of generators of symbolic powers, Waldschmidt constant, (asymptotic) resurgence number are the same as those of group $D_3$, see \cite[Section 4]{ThaiFermat}.\par
\vspace{0.5em}
Geometrically, $J$ is the defining ideal of the singular locus of the line arrangement in $\PP^2$ that consists of $6$ lines $L_j$ whose equations are $$x=0,y=0,z=0,x= y, y= z, z= x$$
These $6$ lines intersect at $7$ points $P_i$ which are $[1:0:0], [0:1:0], [0:0:1], [1:1:1],$ $[0:1:1], [1:0:1]$ and $[1:1:0]$ such that the first $4$ points lie on $3$ lines each and the rest lie on $2$ lines each; and each line contains exactly $3$ points.\par
\vspace{0.5em}
Algebraically, we can write 
$$ J= ((y-z)(y+z-x),(x-y)(x+y-z)) \cap (x,y) \cap (y,z) \cap (z,x)$$ 
and 
$$J^{(m)}=((y-z)(y+z-x),(x-y)(x+y-z))^m \cap (x,y)^m \cap (y,z)^m \cap (z,x)^m \ \ \forall m $$

\begin{theorem}
For ideal $J$ we have the following
\begin{enumerate}
    \item $\widehat{\alpha}(J) = \dfrac{5}{2}$
    \item $\alpha(J^{(2k)})=5k$, for all $k\ge 2$.
    \item $\alpha(J^{(2k+1)})=5k+3$, for all $k\ge 0$.
    \item $\alpha(J^{(2)})=6$
\end{enumerate}
\end{theorem}

\begin{proof}
By \cite[Theorem 2.3]{FaGHLMaS}, we can check that $\widehat{\alpha}(J) \ge \dfrac{5}{2}$.
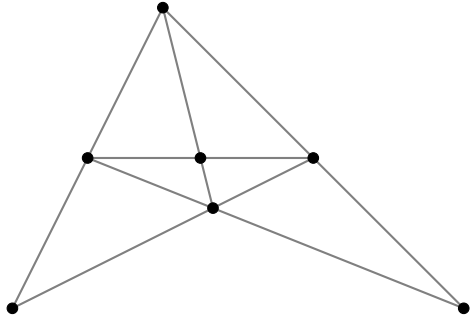
\begin{center}
\begin{tikzpicture}
\draw[gray, thick] (0,4) -- (-2,0);
\draw[gray, thick] (-1,2) -- (2,2);
\draw[gray, thick] (0,4) -- (4,0);
\draw[gray, thick] (0,4) -- (2/3,4/3);
\draw[gray, thick] (4,0) -- (-1,2);
\draw[gray, thick] (-2,0) -- (2,2);
\filldraw[black] (0,4) circle (2pt);
\filldraw[black] (-2,0) circle (2pt);
\filldraw[black] (4,0) circle (2pt);
\filldraw[black] (1/2,2) circle (2pt);
\filldraw[black] (-1,2) circle (2pt);
\filldraw[black] (2,2) circle (2pt);
\filldraw[black] (2/3,4/3) circle (2pt);
\end{tikzpicture}
\end{center}
In particular, we have that $\alpha(J^{(2k)}) \ge 5k$ and $\alpha(J^{(2k+1)}) \ge 5k+3$, for all $k\ge 1$. We will show the reverse by showing there exists some element with the desired degree in the symbolic powers.\par
\vspace{0.5em}
Denote $K=((y-z)(y+z-x),(x-y)(x+y-z))$ and first notice that since $$-(z-x)(z+x-y)=(y-z)(y+z-x)+(x-y)(x+y-z) \in K$$ 
we have that
\begin{equation}\label{inK1_1}
2x(x-y)(z-x)=(z-x)(x-y)(x+y-z)+(x-y)(z-x)(z+x-y) \in K
\end{equation}
Similarly, $y(y-z)(x-y), z(z-x)(y-z)\in K$. It follows that
\begin{equation}\label{inK2}
    (y-z)(z-x)xyz^2=(y-z)(z-x)xy[z^2-(x-y)^2]+(y-z)(z-x)xy(x-y)^2 \in K^2
\end{equation}
since $$(y-z)(z-x)xy[z^2-(x-y)^2]=xy(y-z)(y+z-x)(z-x)(z+x-y) \in K^2$$
and by \ref{inK1_1} $$(y-z)(z-x)xy(x-y)^2=x(x-y)(z-x)y(y-z)(x-y) \in K^2$$
By \ref{inK1_1} we also have that
\begin{equation}\label{inK1_2}
  yz(y-z)=(y-z)z(y+z-x)-(y-z)z(z-x) \in K  
\end{equation}
We have the following cases
\begin{enumerate}
    \item Case 1: when $m=4k$. Consider the polynomial $F=[(y-z)(z-x)xyz^2]^k(x-y)^{2k}(x+y-z)^{2k}$ that has degree $10k$. By \ref{inK1_2}, $[(y-z)(z-x)xyz^2]^k \in K^{2k}$ so $F\in K^{4k}$. On the other hand, $[(y-z)yz^2]^k \in (y,z)^{4k}$, $[(z-x)xz^2]^k \in (z,x)^{4k}$  and $(xy)^{2k}(x-y)^{2k} \in (x,y)^{4k}$. Thus $$F\in K^{4k} \cap (x,y)^{4k} \cap (y,z)^{4k} \cap (z,x)^{4k}=J^{(4k)}, \ \forall k\ge 1$$
    \item Case 2: when $m=4k+2$. We first show that $$G=x^2y^2z^2(x-y)^2(y-z)^2(z-x)^2(x+y-z)(y+z-x)(z+x-y) \in J^{(6)}$$
    In fact, $(x-y)(y-z)(z-x)(x+y-z)(y+z-x)(z+x-y) \in K^3$ and by \ref{inK1_2} $x^2y^2z^2(x-y)^2(y-z)^2(z-x)^2 \in K^3$. It is also clear that $G\in (x,y)^{6} \cap (y,z)^{6} \cap (z,x)^{6} $.\\
    
     The polynomial $F=G[(y-z)(z-x)xyz^2]^k(x-y)^{2k}(x+y-z)^{2k}$ has degree $10(k+1)+5$ and by case $1$ 
     $$F\in J^{(6)}J^{(4k)} \subseteq J^{(4(k+1)+2)}, \ \forall k\ge 0$$
    \item Case 3: when $m=4k+1$. $F=x^ky^{k+1}z^{2k+1}(y-z)^{k+1}(z-x)^{k}(x-y)^{2k}(x+y-z)^{2k}$ has degree $10k+3$. By \ref{inK2}, $x^ky^kz^{2k}(y-z)^k(z-x)^k \in K^{2k}$, by \ref{inK1_2}, $yz(y-z) \in K$ so $F\in K^{4k+1}$. Similar to case $1$, $F \in \cap (x,y)^{4k+1} \cap (y,z)^{4k+1} \cap (z,x)^{4k+1}$ and hence,  
    $$F\in J^{(4k+1)}, \ \forall k\ge 0$$
    \item Case 4: when $m=4k+3$. $F=[(y-z)(z-x)xyz^2]^{k+1}(x-y)^{2k+1}(x+y-z)^{2k+1}$ has degree $10k+8$ and $$F\in K^{4k+3} \cap (x,y)^{4k+3} \cap (y,z)^{4k+3} \cap (z,x)^{4k+3}=J^{(4k+3)}, \ \forall k\ge 0$$
\end{enumerate}

Thus, $\alpha(J^{(2k)})\le 5k$, for all $k\ge 2$ and $\alpha(J^{(2k+1)}) \le 5k+3$, for all $k\ge 0$. It follows that statements $(2)$ and $(3)$ are true and by taking limit as $k$ goes to $\infty$, $(1)$ follows as well. Part $(4)$ can be checked directly by Macaulay2 or by Bezout theorem argument as follows: We know that $\alpha(J^{(2)}) \le \alpha(J^2) =6$. Now suppose that $\alpha(J^{(2)}) \le 5$. Then there is a divisor $D$ of degree $5$ vanishing to order at least $2$ at every point $P_i$. Since the intersection of $D$ and any line $L_j$ consists of $3$ points to order at least $2$, we get a contradiction to Bezout theorem because $\deg(D).\deg(L_j)=5<2.3$.
\end{proof}

\begin{example}
It is worth to point out that the first immediate application of the above calculations is the verification of Chudnovsky's Conjecture and Demailly's Conjecture, although the general case is already known from \cite{EsnaultViehweg}. Ideal $J$ verifies:
\begin{enumerate}
    \item Chudnovsky's Conjecture $\widehat{\alpha}(J)\ge \dfrac{\alpha(J)+1}{2}$.
    \item Demailly's Conjecture $\widehat{\alpha}(J)\ge \dfrac{\alpha(J^{(m)})+1}{m+1}$ for all $m$.
\end{enumerate}
\end{example}

\begin{theorem}
The resurgence number and asymptotic resurgence number are: $$\rho(J)=\widehat{\rho}(J)=\dfrac{6}{5}$$
\end{theorem}
\begin{proof}
Since $J$ has the same minimal degrees of generators of symbolic powers and same free resolutions of all its ordinary powers as those of the ideal of configuration of $D_3$, the proof is identical to that of ideal of configuration of $D_3$ \cite[Proposition 4.2]{ThaiFermat}.
\end{proof}

The above result gives us another example for ideals have expected resurgence. The following corollaries and remark also follow from the fact that ideal of configuration $A_3$ and $D_3$ share the same minimal degrees of generators of all symbolic powers, see \cite[Section 4]{ThaiFermat}.

\begin{corollary}
Ideal $J$ verifies Harbourne-Huneke containment
\begin{enumerate}
    \item $J^{(2r)} \subseteq \mm^rJ^r, \quad \forall r$.
    \item $J^{(2r-1)} \subseteq \mm^{r-1}J^r, \quad \forall r\ge 1$.
\end{enumerate}
 \end{corollary}

\begin{remark}
The above corollary gives a proof for the case $A_3$ in \cite[Proposition 6.3]{DrabkinSeceleanu}.
\end{remark}

In \cite[Example 3.7]{BGHN2020}, we showed the stronger containment (which implies both Harbourne-Huneke) containment $$J^{(2r-2)} \subseteq \mm^{r}J^r$$ for $r=5$ (by Macaulay2) and thus, the same containment holds for all $r \gg 0$ by our method. In particular, from the proof of \cite[Theorem 3.1]{BGHN2020}, the containment hold for $r\ge 10^2=100$. Here we see that the containment hold for all $r\ge 5$. 

\begin{corollary}
For every $n\ge 3$, Fermat configuration ideal verifies the following containment $$J^{(2r-2)} \subseteq \mm^{r}J^r, \quad \forall r\ge 5$$
\end{corollary}

\begin{remark}
For $r\le 4$, the above containment fail with the same reason to that of \cite[Remark 4.6]{ThaiFermat}.
\end{remark}
Table \ref{table:A3} summarizes all numerical values including the initial degrees of symbolic powers as well as the asymptotic resurgence and resurgence numbers of ideal $J$ that we have calculated in this section.
\begin{table}[h!]
  \begin{tabular}{|c||c|c|c|}
    \hline
    $m$ & 2 & $2k$ & $2k+1$\\
    \hline
    $\alpha(J^{(m)})$ & 6 & $5k$ & $5k+3$\\
    \hline
    $\widehat{\alpha}(J)$ &
      \multicolumn{3}{c|}{5/2}\\
    \hline
    $\rho(J) = \widehat{\rho}(J)$ &
      \multicolumn{3}{c|}{6/5}\\
    \hline
  \end{tabular}
  \caption{Least degree of generators and other invariants related to symbolic powers of $J$}
  \label{table:A3}
\end{table}

The table provide a complete answer to the question in the theory of Hermite interpolation, that is to determine the least degree of a homogeneous polynomial that vanishes to order $m$ at the $7$ points of the given configuration in $\PP^2$. \par
\vspace{0.5em}
We end this section by calculating the invariant $\beta(J^{(m)})$. As introduced first in \cite[Definition 2.2]{HaHu} for homogeneous ideals and later considered for ideal $I$ of a finite set of (fat) points in \cite{HaSeFermat}, $\beta(I)$ is set to be the smallest integer $t$ such that $I_t$ contains a regular sequence of length two, or equivalently, is the least degree $t$ such that the zero locus of $I_t$ is $0$-dimensional. It is known that for ideals of fat points, $\beta(I)\le \omega(I)$ \cite[Proposition 3.9]{ThaiFermat}. By the same argument to that of \cite[Proposition 4.9]{ThaiFermat} we have the following

\begin{proposition}
For all $m\ge 1$, $\beta(J^{(m)})=3m$ and $\omega(J^{(m)}) \ge 3m$.
\end{proposition}
\begin{proof}
The proof is identical to that of \cite[Proposition 4.9]{ThaiFermat} since $J^m$ is generated in degree $3m$ and each line in this configuration also passed through exactly $3$ points in the configuration.
\end{proof}

\begin{remark}
Similar to the Fermat ideals, it is suggested by Macaulay2 that $\omega(J^{(m)}) = 3m$.
\end{remark}

\section{An Additional Example}
Consider ideal $J_2=(yz(y^2-z^2),zx(z^2-x^2),xy(x^2-y^2))$. Geometrically, $J$ is the defining ideal of the singular locus of the line arrangement in $\PP^2$ corresponds to the group $G(2,1,3)=B_3$, \cite[Proposition 3.10]{DrabkinSeceleanu}. This arrangement consists of $9$ lines $L_j$ whose equations are $$x=0,y=0,z=0,x= \pm y, y= \pm z, z= \pm x$$
These $9$ lines intersect at $13$ points $P_i$ which are $[1:0:0], [0:1:0], [0:0:1], [1:1:1],$ $[-1:1:1], [1:-1:1], [1:1:-1]$, $[0:1:1],[1:0:1],[1:1:0]$, $[0:1:-1],[-1:0:1]$ and $[1:-1:0]$ such that the first $3$ points lie on $4$ lines each, the next $4$ points lie on $3$ lines each and the rest lie on $2$ lines each; and each line contains exactly $4$ points.\par
\vspace{0.5em}
Algebraically, we can write 
$$ J_2= K_1 \cap K_2 \cap K_3 \cap K \cap (x,y) \cap (y,z) \cap (z,x)$$
where $K_1=(x,y^2-z^2),K_2=(y,z^2-x^2),K_3=(z,x^2-y^2)$ and $K=(x^2-y^2,y^2-z^2)$ ,or we can write 
$$ J_2= K'_1 \cap K'_2 \cap K \cap (x,y) \cap (y,z) \cap (z,x)$$
where $K'_1=((y-z)(y+z-x),(x-y)(x+y-z)), K'_2=(x+y+z,yz(y+z))$, then for all $m$
\begin{align*}
    J_2^{(m)} &= K_1^m \cap K_2^m \cap K_3^m \cap K^m \cap (x,y)^m \cap (y,z)^m \cap (z,x)^m \\
    &= K_1^{'m} \cap K_2^{'m} \cap K^m \cap (x,y)^m \cap (y,z)^m \cap (z,x)^m \\
\end{align*}

Note that $K \cap (x,y) \cap (y,z) \cap (z,x)$ is the Fermat ideal with $n=2$, so we immediately have that $\widehat{\alpha}(J_2) \ge \dfrac{5}{2}$. Moreover, computations with Macaulay2 show that $\alpha(J_2^{(6)})=21$, hence, $\widehat{\alpha}(J_2) \le \dfrac{7}{2}$. It seems that $\widehat{\alpha}(J_2)$ would be $\dfrac{7}{2}$ but the $\alpha(J_2^{(m)})$ are tricky to deal with. On the other hand, $\alpha(J_2^{(3)})=12$, by \cite{EsnaultViehweg}, since $J_2$ satisfies Demailly's Conjecture, we have: 
$$\widehat{\alpha}(J_2)\ge \dfrac{\alpha(J_2^{(3)})+1}{3+1}=\dfrac{13}{4}$$
Since $\reg(J_2)=6$, by \cite[Theorem 1.2]{AsymptoticResurgence} we have:
$$\dfrac{8}{7} \le \dfrac{\alpha(J_2)}{\widehat{\alpha}(J_2)} \le \widehat{\rho}(J_2) \le \rho(J_2) \le \dfrac{\reg(J_2)}{\widehat{\alpha}(J_2)} \le \dfrac{24}{13}$$

In particular, $J_2$ has expected resurgence. Further computations by Macaulay2 show that, $J_2^{(7)} \not \subset J_2^6$, thus, $\dfrac{7}{6}\le \rho(J_2)$. It is interesting to know if $\rho(J_2) = \dfrac{7}{6}$.\par
\vspace{0.5em}
More general, for $n\ge 1$, consider the ideal $J_n=(yz(y^n-z^n),zx(z^n-x^n),xy(x^n-y^n))$ which will capture $J_2$ and $J=J_1$. Geometrically, $J_n$ corresponds to the configuration of $3n+3$ lines that consists of all $3n$ lines in Fermat configuration and $3$ lines $x=0,y=0,z=0$; and $n^2+3n+3$ intersection points. In particular, each line passes through $n+2$ points. Similar to $J_2$, note that $J_n$ is a subset of the Fermat ideals for $n\ge 2$, thus $\widehat{\alpha}(J_n) \ge n$ for $n\ge 3$. 
\begin{remark}
Only by looking at the above lower bound, we easily see that $J_n$ verifies Chudnovsky's Conjecture $\widehat{\alpha}(J_n)\ge \dfrac{\alpha(J_n)+1}{2}$ for all $n$. In fact, the case $n=1$ was verified in previous sections and $n=2$ follows from $\widehat{\alpha}(J_2) \ge \dfrac{5}{2}$. When $n \ge 3$, $\widehat{\alpha}(J_n) \ge n \ge \dfrac{n+2+1}{2}$.
\end{remark}

On the other hand, by \cite[Theorem 2.5]{NagelSeceleanu}, since $J_n$ is a strict almost complete intersection ideal with minimal generators of degree $n+2$ and its module sygyzies is generated in degree $1$ and $n+1$, the minimal free resolution of $J_n^r$ is:
$$ 0 \to R(-(n+2)(r+1))^{{r} \choose 2}\stackrel{\psi}\to \begin{matrix} R(-(n+2)r-1)^{{r+1} \choose 2} \\
\oplus \\ R(-(n+2)r-(n+1))^{{r+1} \choose 2}\end{matrix} \stackrel{\varphi}\to R(-(n+2)r)^{{r+2} \choose 2} \to J_n^r\to 0 $$
for any $r\ge 2$, in particular, $\reg(J_n^r)=(n+2)r+n$ for all $r\ge 2$. It would be interesting to know if we can determine $\alpha(J_n^{(m)})$ and use them with the knowledge of $\reg(J_n^r)$ to verify Demailly's Conjecture as well as Harbourne-Huneke Containment, stable Harbourne Containment as we did for Fermat ideals and Fermat-like ideals.\par
\vspace{0.5em}
Back to ideal $J_2$, the following proposition (part $(2)$) gives a proof for the case $B_3$ in \cite[Proposition 6.3]{DrabkinSeceleanu}.

\begin{proposition}
Ideal $J_2$ verifies Harbourne-Huneke containment
\begin{enumerate}
    \item $J_2^{(2r)} \subseteq \mm^rJ_2^r, \quad \forall r$.
    \item $J_2^{(2r-1)} \subseteq \mm^{r-1}J_2^r, \quad \forall r\ge 3$.
\end{enumerate}
\end{proposition}
\begin{proof}
\begin{enumerate}
    \item For all $r$, $J_2^{(2r)} \subseteq J_2^r$, hence, the containment follows since for all $r$, we have
    $$\alpha(J_2^{(2r)}) \ge \dfrac{13r}{2} \ge r+4r = r+\omega(J_2^r)$$
    \item First, the containment $J_2^{(2r-1)} \subseteq J_2^r$ for $r\ge 3$ follows from the inequality 
    $$\alpha(J_2^{(2r-1)}) \ge\dfrac{13}{4}(2r-1)\ge 4r+2 = \reg(J_2^r)$$
    for all $r\ge 3$. Thus, $J_2^{(2r-1)} \subseteq \mm^{r-1}J_2^r$ for $r\ge 3$ follows from: 
    $$\alpha(J_2^{(2r-1)}) \ge\dfrac{13}{4}(2r-1)\ge r-1+4r = r+\omega(J_2^r)$$
    for all $r\ge 2$.
\end{enumerate}
\end{proof}

\begin{remark}
Macaulay2 shows that for $r=2$, $\alpha(J_2^{(3)}) =12$, hence, $\alpha(J_2^{(3)}) \ge 4.2+2 = \reg(J_2^2)$. Therefore, $J_2^{(2r-1)} \subseteq J_2^r$ for $r\ge 2$ and since $\alpha(J_2^{(2r-1)}) \ge r+\omega(J_2^r)$ for all $r\ge 2$, $J_2^{(2r-1)} \subseteq \mm^{r-1}J_2^r$, for $r\ge 2$. The case $r=1$ is obvious. Thus, $J_2^{(2r-1)} \subseteq \mm^{r-1}J_2^r$, for $r\ge 1$.
\end{remark}

We end this section with the follow up to the discussion about the invariants $\beta$ and $\omega$ for Fermat ideals in \cite{ThaiFermat} and $J=J_1$ in the previous section.

\begin{proposition}
For all $n\ge 1$ and $m\ge 1$, $\beta(J_n^{(m)})=m(n+2)$ and $\omega(J_n^{(m)}) \ge m(n+2)$.
\end{proposition}
\begin{proof}
The proof is the same to that of \cite[Proposition 3.10]{ThaiFermat} with notice that $J_n^m$ is generated in degree $m(n+2)$ and each line in this configuration passed through exactly $n+2$ points in the configuration.
\end{proof}

\begin{remark}
It is also suggested by Macaulay2 \cite{M2} that $\omega(J_n^{(m)}) = m(n+2)$ as in the case of Fermat ideals. It is interesting to know if $\omega(I^{(m)})=\beta(I^{(m)})$ hold for what radical ideal of points $I$ in general.
\end{remark}


\end{document}